# Spectral integrated neural networks (SINNs) for solving forward and inverse dynamic problems

Lin Qiu [a], Fajie Wang [a], Wenzhen Qu [b,*], Yan Gu [b,*], Qing-Hua Qin [c]

[a] College of Mechanical and Electrical Engineering, National Engineering Research Center for Intelligent Electrical Vehicle Power System, Qingdao University, Qingdao 266071, PR China

[b] School of Mathematics and Statistics, Qingdao University, Qingdao 266071, PR China

[c] Department of Materials Science, Shenzhen MSU-BIT University, Shenzhen 518172, PR China

**Abstract**

This paper proposes a novel neural network framework, denoted as spectral integrated neural networks (SINNs), for resolving three-dimensional forward and inverse dynamic problems. In the SINNs, the spectral integration method is applied to perform temporal discretization, and then a fully connected neural network is adopted to solve resulting partial differential equations (PDEs) in the spatial domain. Specifically, spatial coordinates are employed as inputs in the network architecture, and the output layer is configured with multiple outputs, each dedicated to approximating solutions at different time instances characterized by Gaussian points used in the spectral method. By leveraging the automatic differentiation technique and spectral integration scheme, the SINNs minimize the loss function, constructed based on the governing PDEs and boundary conditions, to obtain solutions for dynamic problems. Additionally, we utilize polynomial basis functions to expand the unknown function, aiming to enhance the performance of SINNs in addressing inverse problems. The conceived framework is tested on six forward and inverse dynamic problems, involving nonlinear PDEs. Numerical results demonstrate the superior performance of SINNs over the popularly used physics-informed neural networks in terms of convergence speed, computational accuracy and efficiency. It is also noteworthy that the SINNs exhibit the capability to deliver accurate and stable solutions for long-time dynamic problems.

---
[*] Corresponding authors.
 *E-mail addresses:* quwz@qdu.edu.cn (W. Qu); guyan1913@163.com (Y. Gu).

# 1. Introduction

Numerical simulations play a pivotal role in engineering and scientific research, providing insight into complex physical phenomena. Some numerical methodologies, such as the finite element method (FEM) [1-3], boundary element method (BEM) [4], and meshless methods [5-8] like smoothed particle hydrodynamics (SPH) [9], element-free Galerkin (EFG) [10], etc., have emerged as effective tools for simulating diverse phenomena. These approaches exhibit individual advantages in numerical simulations. However, inherent limitations or disadvantages persist in these techniques that need further addressing, such as the troublesome mesh generation in mesh-based schemes, the intricate singular integrals in the BEM, and the poor stability in some meshless methods.

In recent years, propelled by advancements in computational resources and optimization approaches, a novel machine learning architecture named physics-informed neural networks (PINNs) [11, 12] has become a focal point of research interest. The PINNs leverage the advantages of both deep learning and traditional modeling methods by embedding the underlying partial differential equations (PDEs), i.e. physical laws of the problem, into the structure of neural networks. This fusion enables PINNs to be a powerful tool for tackling intricate problems in the realm of mechanics [13-15], engineering [16], materials science [17], biomedicine [18], energy [19, 20], etc. Various variants of the PINNs are conceived to facilitate diverse applications, including conservative PINNs [21], fractional PINNs [22], Bayesian PINNs [23], and variational PINNs [24]. To enhance the performance of PINNs, some tests on different activation functions [25, 26], network architectures [27, 28], and optimization algorithms [29, 30] are conducted. Additionally, some strategies are proposed and confirmed to be effective, such as domain decomposition [31], self-adaptive approach [32], importance sampling [33], etc. Despite significant progress has been achieved in the theory and applications of PINNs, their convergence, stability, computational efficiency and accuracy still to some extent restrict the applications in practical scenarios, especially in dealing with high-dimensional dynamic problems.

Currently, the PINNs typically take spatiotemporal information as network input and employ automatic differentiation technique to directly handle the temporal derivative for solving time-dependent problems [34]. Alternatively, low-order differentiation scheme can be combined to discretize the temporal term. These approaches are found to incur significant computational costs and exhibit limited computational accuracy. In particular, it is extremely challenging for such schemes to stably perform long-time numerical simulations. The spectral integration method [35], relying on Gaussian quadrature and the basis functions like Legendre polynomial and Chebyshev polynomial, offers a potential strategy to address this challenge, due to its rapid convergence, high accuracy and stability. With these advantages, this method has been successfully applied to assist in solving many types of PDEs [35, 36]. Moreover, this method has spawned various approaches, including the spectral deferred correction (SDC) method [37] and Krylov deferred correction (KDC) method [38], and has achieved extensive applications.

Leveraging the principles of PINNs and spectral integration method, we establish a novel neural network framework, referred as spectral integrated neural networks (SINNs), for accurately and efficiently dealing with three-dimensional (3D) forward and inverse dynamic problems. In this approach, the spectral integration method is utilized for temporal discretization of dynamic problems, and a fully connected neural network is employed to obtain solutions in the spatial domain. Specifically, the SINNs employ spatial coordinates as network inputs, and configure multiple outputs in the output layer to approximate solutions at different time instances identified by Gaussian points used in the spectral integration method. The outputs are processed by adopting the automatic differentiation technique and spectral integration scheme, and then a loss function is constructed based on the governing PDEs and boundary conditions. The parameters of the neural network are trained by using the back-propagation of loss function and the gradient descent method. To improve the efficacy of SINNs in tackling inverse dynamic problems, we expand the function to be recovered by utilizing polynomial basis functions with unknown coefficients. By embedding these coefficients into the architecture of the neural network, the unknown coefficients and network parameters are updated and trained simultaneously. The performance of the SINNs is evaluated by several representative dynamic problems, encompassing linear and nonlinear transient heat conduction in functionally graded materials (FGM) [39, 40], linear and nonlinear transient wave propagation problems, inverse problem of heat conduction in FGM, and long-time heat conduction in FGM. Furthermore, we conduct a comparison between the PINNs and proposed SINNs with an emphasis on convergence speed, accuracy, and efficiency.

The rest of the paper is structured as follows. Section 2 introduces the 3D dynamic problems, focusing on the mathematical models of transient heat conduction and wave propagation problems. In Section 3, we present the spectral integration method, and give the formulation and schematic of SINNs for resolving the forward and inverse dynamic problems. In Section 4, the performance of the proposed SINNs is examined through six benchmark examples related to the linear and nonlinear dynamic problems. Finally, some conclusions and remarks are summarized in Section 5.

**2. Problem setup**

Assuming that $\Omega$ is a 3D space domain bounded by a surface $\partial\Omega = \Gamma$, where $\Gamma = \Gamma_D \cup \Gamma_N$ and $\Gamma_D \cap \Gamma_N = \varnothing$, we define the general form of dynamic problems as follows,

$$\mathcal{L}_t\big[u(\boldsymbol{x},t);\lambda(\cdot)\big] + \mathcal{L}_s\big[u(\boldsymbol{x},t);\mu(\cdot)\big] = f(\boldsymbol{x},t), \ \boldsymbol{x} \in \Omega, t \in [0,T], \tag{1}$$

where $\boldsymbol{x} = x,y,z$ is the spatial coordinate, $u(\boldsymbol{x},t)$ stands for the latent solution, $f(\boldsymbol{x},t)$ indicates the source term, $\mathcal{L}_t$ and $\mathcal{L}_s$ denote general differential operators that are constructed with time and space derivatives, respectively. $\lambda(\cdot)$ and $\mu(\cdot)$ represent parameters, which may be either constants or functions associated with $\boldsymbol{x}$ and $t$. In this work, we will focus on two types of dynamics problems, i.e. the transient heat conduction and wave propagation problems.

For transient heat conduction problems, the temperature function $u(\boldsymbol{x},t)$ satisfy the following governing equation [41, 42],

$$\rho(\cdot)c(\cdot)u_t(\boldsymbol{x},t) - \nabla\big[\kappa(\cdot)\nabla u(\boldsymbol{x},t)\big] = f(\boldsymbol{x},t), \ \boldsymbol{x} \in \Omega, t \in [0,T], \tag{2}$$

subject to the initial condition, Dirichlet and Neumann boundary conditions,

$$\begin{cases} u(\boldsymbol{x},0) = u_0(\boldsymbol{x}), \ \boldsymbol{x} \in \Omega, \\ u(\boldsymbol{x},t) = \overline{u}(\boldsymbol{x},t), \ \boldsymbol{x} \in \Gamma_D, t \in [0,T], \\ q(\boldsymbol{x},t) = -\kappa(\cdot)\dfrac{\partial u(\boldsymbol{x},t)}{\partial \boldsymbol{n}} = \overline{q}(\boldsymbol{x},t), \ \boldsymbol{x} \in \Gamma_N, t \in [0,T], \end{cases} \tag{3}$$

where $\rho(\cdot)$, $c(\cdot)$ and $\kappa(\cdot)$ denote the mass density, the specific heat and the thermal conductivity, respectively. $f(\boldsymbol{x},t)$ represents the source term. In particular, Eq. (2) characterizes the heat transfer in FGM when $\rho(\cdot)$, $c(\cdot)$ and $\kappa(\cdot)$ are all functions related to $\boldsymbol{x}$, and it describes nonlinear transient heat conduction problems when $\kappa(\cdot)$ is a function of the temperature function $u(\boldsymbol{x},t)$ [43].

For transient wave propagation problems, the physical variable $u(\boldsymbol{x},t)$ satisfy the following governing equation [44, 45],

$$u_{tt}(\boldsymbol{x},t) - w^2\nabla^2 u(\boldsymbol{x},t) + \mathcal{N}\big[u(\boldsymbol{x},t);\hat{\mu}(\cdot)\big] = f(\boldsymbol{x},t), \ \boldsymbol{x} \in \Omega, t \in [0,T], \tag{4}$$

subject to the initial conditions, Dirichlet and Neumann boundary conditions,

$$\begin{cases} u(\boldsymbol{x},0) = u_0(\boldsymbol{x}), \ \boldsymbol{x} \in \Omega, \\ u_t(\boldsymbol{x},0) = v_0(\boldsymbol{x}), \ \boldsymbol{x} \in \Omega, \\ u(\boldsymbol{x},t) = \overline{u}(\boldsymbol{x},t), \ \boldsymbol{x} \in \Gamma_{\mathrm{D}}, t \in [0,T], \\ q(\boldsymbol{x},t) = \dfrac{\partial u(\boldsymbol{x},t)}{\partial \boldsymbol{n}} = \overline{q}(\boldsymbol{x},t), \ \boldsymbol{x} \in \Gamma_{\mathrm{N}}, t \in [0,T], \end{cases} \quad (5)$$

where $w$ denotes the wave propagation speed, $\mathcal{N}\big[u(\boldsymbol{x},t);\hat{\mu}(\cdot)\big]$ represents a nonlinear term parameterized by $\hat{\mu}(\cdot)$.

## 3. Methodology

### 3.1. Spectral integration method

Following the idea of the Picard integral equation commonly employed in spectral methods [37, 38], we introduce $U(\boldsymbol{x},t) = u_t(\boldsymbol{x},t)$ as the new variable to be determined instead of directly resolving $u(\boldsymbol{x},t)$ in Eq. (2). The transient heat conduction equation is reformulated as

$$\rho(\cdot)c(\cdot)U(\boldsymbol{x},t) - \nabla \cdot \kappa(\cdot)\nabla\left[u_0(\boldsymbol{x}) + \int_0^t U(\boldsymbol{x},\tau)d\tau\right] = f(\boldsymbol{x},t), \ \boldsymbol{x} \in \Omega, t \in [0,T], \quad (6)$$

where $u_0(\boldsymbol{x}) = u(\boldsymbol{x},0)$ indicates the initial value. When encountering a problem with a long time interval that is difficult to solve directly, we divide that interval into multiple subintervals, i.e. $[0,T] = [T_0 = 0, T_1] \cup [T_1, T_2] \cup \cdots \cup [T_{n-1}, T_n = T]$, and denote the step size as $\Delta T_i = T_i - T_{i-1}$. The key lies in effectively resolving the problem within the small interval $[0, T_1]$. Subsequently, leveraging the time-marching scheme [46] and using the values at $T_1$, the subsequent subintervals can be addressed sequentially.

More generally, we clarify how to solve Eq. (6) within the subinterval $[T_{i-1}, T_i]$, which can be mapped to $[-1,1]$ by a linear transformation. We arrange $p$ Gaussian-type nodes $\{\tilde{T}_j\}_{j=1}^p$ in the subinterval, and mark the corresponding values of the new variable at the nodes as $\boldsymbol{U} = \big[U(\boldsymbol{x},\tilde{T}_1), U(\boldsymbol{x},\tilde{T}_2), \cdots, U(\boldsymbol{x},\tilde{T}_p)\big]$. The value $\boldsymbol{U}$ can be stably and accurately approximated by constructing Legendre polynomial expansion $L^p(\boldsymbol{U},t)$, the coefficients of which are determined utilizing Gaussian quadrature rules [38]. We integrate the interpolation polynomial from $T_{i-1}$ to $\tilde{T}_j$, and further derive the mapping relationship between the value $\boldsymbol{U}$ and the integral as follows,

$$\left[\Delta T_i S \boldsymbol{U}\right]_j = \int_{T_{i-1}}^{\tilde{T}_j} L^p(\boldsymbol{U},\tau)d\tau, \tag{7}$$

where $S$ stands for the spectral integration matrix independent of $\Delta T_i$. For a more detailed description regarding the spectral integration matrix, interested readers are advised to refer to Refs. [35, 36, 38]. Utilizing the spectral integration matrix, the governing equation is discretized in time as

$$\rho(\cdot)c(\cdot)U(\boldsymbol{x},\tilde{T}_j) - \nabla\left\langle \kappa(\cdot)\nabla\ u(\boldsymbol{x},T_{i-1}) + \left[\Delta T_i S\boldsymbol{U}\right]_j\ \right\rangle = f(\boldsymbol{x},\tilde{T}_j), \boldsymbol{x}\in\Omega, j=1,2,\cdots,p, \tag{8}$$

and the corresponding boundary conditions in Eq. (3) are converted into the following expression,

$$\begin{cases} U(\boldsymbol{x},\tilde{T}_j) = \bar{u}_t(\boldsymbol{x},\tilde{T}_j) = \bar{U}(\boldsymbol{x},\tilde{T}_j),\ \boldsymbol{x}\in\Gamma_{\mathrm{D}}, j=1,2,\cdots,p, \\ Q(\boldsymbol{x},\tilde{T}_j) = -\kappa(\cdot)\dfrac{\partial U(\boldsymbol{x},\tilde{T}_j)}{\partial \boldsymbol{n}} = \bar{q}_t(\boldsymbol{x},\tilde{T}_j) = \bar{Q}(\boldsymbol{x},\tilde{T}_j),\ \boldsymbol{x}\in\Gamma_{\mathrm{N}}, j=1,2,\cdots,p, \end{cases} \tag{9}$$

where $f(\boldsymbol{x},\tilde{T}_j)$ represents the value of source term at $j$ th Gaussian node, $\bar{U}(\boldsymbol{x},\tilde{T}_j)$ and $\bar{Q}(\boldsymbol{x},\tilde{T}_j)$ indicate known values at $j$ th Gaussian node.

For transient wave propagation problems involving second-order time derivatives, we define $U(\boldsymbol{x},t) = u_{tt}(\boldsymbol{x},t)$ as the new variable, and we have,

$$u(\boldsymbol{x},t) = u_0(\boldsymbol{x}) + v_0(\boldsymbol{x})t + \int_0^t \int_0^\tau U(\boldsymbol{x},\varsigma)d\varsigma d\tau, \tag{10}$$

in which $u_0(\boldsymbol{x}) = u(\boldsymbol{x},0)$ and $v_0(\boldsymbol{x}) = u_t(\boldsymbol{x},0)$ are the initial values. Substituting Eq. (10) into Eq. (4) and drawing an analogy to the treatment of heat conduction problems, the transient wave propagation equation is discretized in time with Gaussian nodes, and the boundary conditions can be derived.

*3.2. Formulation and architecture of spectral integrated neural networks*
*3.2.1. Spectral integrated neural networks for forward dynamic problems*

In this section, we illustrate the formulation and framework of SINNs for resolving forward dynamic problems, using the transient heat conduction problem as an example. The schematic of SINNs architecture is illustrated in Fig. 1. Similar to the deep learning architecture outlined in Ref. [11], the developed SINNs numerically approximate the solution of the problem utilizing a fully connected network. The primary distinction is that, when addressing dynamic problems, the SINNs do not incorporate the time coordinate as an initial input in the input layer. Instead, only spatial coordinates are employed as inputs, and multiple outputs are configured in the output layer to represent the values of the trial solution at different Gaussian nodes. During the process of data transmission in the neural network, it is adjusted by weights $\boldsymbol{w}$ and biases $\boldsymbol{b}$, as well as the activation functions $\sigma$. Different outputs of SINNs are assigned distinct network parameters, which are individually labeled as $\boldsymbol{\theta}^{(j)} = \{\boldsymbol{w}^{(j)}, \boldsymbol{b}^{(j)}\}, j = 1, 2, \cdots, p$. The corresponding trial solutions are denoted as $U(\boldsymbol{x}, \tilde{T}_j; \boldsymbol{\theta}^{(j)}), j = 1, 2, \cdots, p$.

With the view of using Eq. (8), we define the PDE residuals as,

$$\mathcal{R}_{PDE}^{(j)}\left(\boldsymbol{x}; \boldsymbol{\theta}^{(1)-(p)}\right) := \nabla \left\langle \kappa(\cdot) \nabla \left\{ u(\boldsymbol{x}, T_{i-1}) + \left[\Delta T_i S \boldsymbol{U}(\boldsymbol{\theta}^{(1)-(p)})\right]_j \right\} \right\rangle \\ + f(\boldsymbol{x}, \tilde{T}_j) - \rho(\cdot)c(\cdot)U(\boldsymbol{x}, \tilde{T}_j; \boldsymbol{\theta}^{(j)}), j = 1, 2, \cdots, p, \quad (11)$$

where $\boldsymbol{\theta}^{(1)-(p)} = (\boldsymbol{\theta}^{(1)}; \boldsymbol{\theta}^{(2)}; \cdots; \boldsymbol{\theta}^{(p)})$. It is noted that only spatial derivatives are involved in the SINNs, and all spatial derivatives can be easily calculated by adopting the automatic differentiation technique. In the SINNs, the governing PDE and boundary conditions are introduced to constrain the trial solutions for training the parameters of PINNs. This idea is implemented by minimizing the following loss function,

$$Loss\left(\boldsymbol{\theta}^{(1)-(p)}; N\right) = \sum_{j=1}^{p} Loss_{PDE}^{(j)}\left(\boldsymbol{\theta}^{(1)-(p)}; N_{PDE}\right) + \sum_{j=1}^{p} Loss_{DBC}^{(j)}\left(\boldsymbol{\theta}^{(j)}; N_{DBC}\right) \\ + \sum_{j=1}^{p} Loss_{NBC}^{(j)}\left(\boldsymbol{\theta}^{(j)}; N_{NBC}\right), \quad (12)$$

in which $Loss_{PDE}^{(j)}$, $Loss_{DBC}^{(j)}$ and $Loss_{NBC}^{(j)}$ stand for the losses of governing PDE, Dirichlet and Neumann boundary conditions at $j$ th Gaussian node, respectively. $N_{PDE}$, $N_{NBC}$ and $N_{DBC}$ represent the number of nodes involved in the corresponding items. The loss terms are detailed below,

$$Loss_{PDE}^{(j)}\left(\boldsymbol{\theta}^{(1)-(p)}; N_{PDE}\right) = \frac{1}{N_{PDE}} \sum_{k=1}^{N_{PDE}} \left|\mathcal{R}_{PDE}^{(j)}\left(\boldsymbol{x}_{PDE}^k; \boldsymbol{\theta}^{(1)-(p)}\right)\right|^2, j = 1, 2, \cdots, p, \quad (13)$$

$$Loss_{DBC}^{(j)}\left(\boldsymbol{\theta}^{(j)}; N_{DBC}\right) = \frac{1}{N_{DBC}} \sum_{k=1}^{N_{DBC}} \left|U\left(\boldsymbol{x}_{DBC}^k, \tilde{T}_j; \boldsymbol{\theta}^{(j)}\right) - \bar{U}\left(\boldsymbol{x}_{DBC}^k, \tilde{T}_j\right)\right|^2, j = 1, 2, \cdots, p. \quad (14)$$

$$Loss_{NBC}^{(j)}\left(\boldsymbol{\theta}^{(j)}; N_{NBC}\right) = \frac{1}{N_{NBC}} \sum_{k=1}^{N_{NBC}} \left|\kappa(\cdot)\frac{\partial U\left(\boldsymbol{x}_{NBC}^k, \tilde{T}_j; \boldsymbol{\theta}^{(j)}\right)}{\partial \boldsymbol{n}} + \bar{Q}\left(\boldsymbol{x}^k, \tilde{T}_j\right)\right|^2, j = 1, 2, \cdots, p, \quad (15)$$

where $\{\boldsymbol{x}_{PDE}^k\}_{k=1}^{N_{PDE}}$ denote collocation points governed by the PDE, $\{\boldsymbol{x}_{DBC}^k\}_{k=1}^{N_{DBC}}$ and $\{\boldsymbol{x}_{NBC}^k\}_{k=1}^{N_{NBC}}$ indicate training points subject to Dirichlet and Neumann boundary conditions, respectively. The parameters of the fully connected network are learned by utilizing the gradient descent method based on the back-propagation of loss function. Subsequently, the trained network is employed to predict the values of physical variable throughout the computational domain.

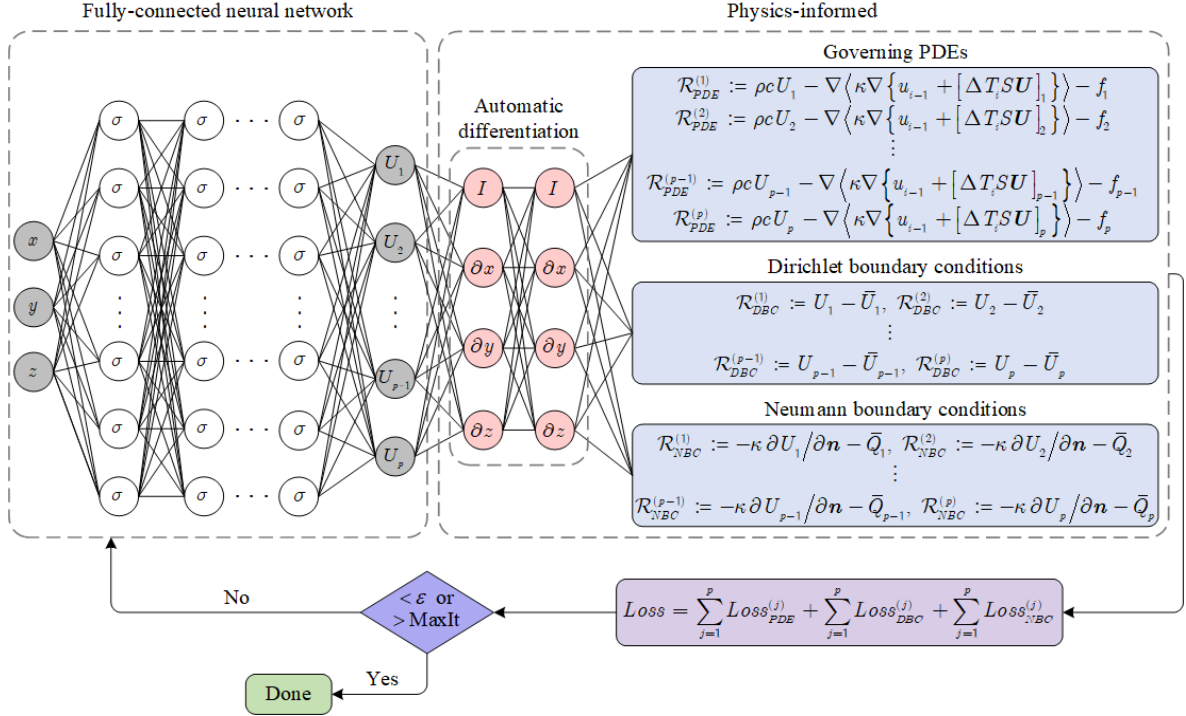

**Fig. 1.** Schematic diagram of the SINNs.

*3.2.2. Spectral integrated neural networks for inverse dynamic problems*

In this section, we aim to employ the proposed SINNs to simultaneously recover the physical variable and the material parameters, such as thermal conductivity, the mass density and the specific heat, portrayed by a function $d(\boldsymbol{x})$. Generally, some data needs to be overspecified to address inverse problems. In this study, additional boundary data are measured and utilized. To enhance the performance of the SINNs in solving inverse problems, this paper adopts the following combination form of polynomial basis functions with order $s$ to represent the function to be inverted,

$$d(\boldsymbol{x}) = \sum_{p=0}^{s}\sum_{q=0}^{p}\sum_{r=0}^{p-q}\alpha_{pqr}\mathcal{B}_{pqr}^{s}(\boldsymbol{x}), \tag{16}$$

in which $\{\alpha_{pqr}\}$ are the unknown coefficients, marked as $\boldsymbol{\alpha}$, and

$$\mathcal{B}_{pqr}^{s}(\boldsymbol{x}) = \{x^{p-q-r}y^{q}z^{r} : 0 \leq p \leq s, 0 \leq q \leq p, 0 \leq r \leq p-q\}. \tag{17}$$

It is noted that the number of polynomial basis functions is $\frac{1}{6}(s+1)(s+2)(s+3)$. By embedding the coefficients $\boldsymbol{\alpha}$ into the structure of the neural network, the unknown coefficients and network parameters are synchronously updated and trained. This process restricts the loss function to gradually approach zero, thereby simultaneously obtaining both physical variables and material parameters. For problems involving multiple material parameters related to the function $d(\boldsymbol{x})$, such as the heat conduction in FGM, the material parameters take the following form,

$$\kappa(\boldsymbol{x}) = \kappa_0 d(\boldsymbol{x}), \rho(\boldsymbol{x})c(\boldsymbol{x}) = \rho_0 c_0 d(\boldsymbol{x}), \tag{18}$$

where $\kappa_0$ and $\rho_0 c_0$ are constants, we can incorporate different coefficients, $\lambda_1$ and $\lambda_2$, into the network to assist in approximating $\kappa(\boldsymbol{x})$ and $\rho(\boldsymbol{x})c(\boldsymbol{x})$, respectively.

In SINNs, the formulation established for the forward problems can be almost directly applied to solve the corresponding inverse problems. The loss function $Loss_{inv}$ is defined as,

$$\begin{aligned} Loss_{inv}\left(\boldsymbol{\theta}^{(1)-(p)};\boldsymbol{\alpha};\lambda_1,\lambda_2;N\right) &= \sum_{j=1}^{p} Loss_{PDE,inv}^{(j)}\left(\boldsymbol{\theta}^{(1)-(p)};\boldsymbol{\alpha};\lambda_1,\lambda_2;N_{PDE}\right) \\ &+ \sum_{j=1}^{p} Loss_{DBC,inv}^{(j)}\left(\boldsymbol{\theta}^{(j)};N_{DBC}\right) + \sum_{j=1}^{p} Loss_{NBC,inv}^{(j)}\left(\boldsymbol{\theta}^{(j)};\boldsymbol{\alpha};\lambda_1;N_{NBC}\right). \end{aligned} \tag{19}$$

These loss terms in Eq. (19) are similar to those described in forward problems. The difference lies in that, compared to Eq. (12), Eq. (19) will contain more boundary data. Specifically, in this paper, Dirichlet boundary conditions on all boundaries or surfaces and Neumann boundary conditions on some boundaries or surfaces are known.

**4. Numerical examples and discussions**

To assess the performance of the proposed SINNs in dealing with both forward and inverse dynamic problems, we test the SINNs on several numerical examples associated with dynamic heat conduction and wave propagation problems. We compare the numerical results achieved by utilizing the SINNs with those obtained by the widely used PINNs, with a focus on the accuracy and efficiency of the algorithm. All computations were conducted using MATLAB R2022b on a Windows 10 (64-bit) platform, equipped with an Intel Core i7-13700K 3.40 GHz CPU and 128 GB RAM. To measure the numerical accuracy of the method, the following relative error and $L_2$ relative error are defined,

$$\text{Relative error} = \left|u(\boldsymbol{x}_m,t_m) - \tilde{u}(\boldsymbol{x}_m,t_m)\right| / \left|u(\boldsymbol{x}_m,t_m)\right|, \tag{20}$$

$$L_2 \text{ relative error} = \left\{\sum_{m=1}^{N_{test}}\left[u(\boldsymbol{x}_m,t_m) - \tilde{u}(\boldsymbol{x}_m,t_m)\right]^2\right\}^{1/2} \bigg/ \left\{\sum_{m=1}^{N_{test}}\left[u(\boldsymbol{x}_m,t_m)\right]^2\right\}^{1/2}, \quad (21)$$

where $(\boldsymbol{x}_m, t_m)$ is the $m$ th test node, $u(\boldsymbol{x}_m, t_m)$ and $\tilde{u}(\boldsymbol{x}_m, t_m)$ respectively stand for fabricated/exact and numerical solutions, and $N_{test}$ indicates the total number of test nodes.

*4.1. Heat conduction problem in FGM*

We first consider a dynamic heat conduction problem in the gear-shaped FGM, the principal dimension of the gear is $0.60\text{m} \times 0.60\text{m} \times 0.10\text{m}$, as presented in Figs. 2(a)-2(c). The gear-shaped FGM is constrained by mixed-type boundary conditions, where the heat fluxes are enforced on the left surface of the gear ($0\text{m} \leq x \leq 0.25\text{m}$), and the temperatures are prescribed on the remaining surface of the gear. The governing equation is described as,

$$\rho(\boldsymbol{x})c(\boldsymbol{x})u_t(\boldsymbol{x},t) - \nabla\left[\kappa(\boldsymbol{x})\nabla u(\boldsymbol{x},t)\right] = f(\boldsymbol{x},t), \quad \boldsymbol{x} \in \Omega, t \in [0\text{s},1\text{s}], \quad (22)$$

where the material parameters $\rho(\boldsymbol{x})c(\boldsymbol{x}) = 2.2d(\boldsymbol{x})$, $\kappa(\boldsymbol{x}) = 1.3d(\boldsymbol{x})$, and

$$d(\boldsymbol{x}) = 0.5\cos(2x) + 0.2\sin(3y) + 0.2\cos(z) + 1. \quad (23)$$

The fabricated solution for this problem is

$$u(\boldsymbol{x},t) = 30\left[\sin(t) + 0.5\right](x + 2y + 3z)^2, \quad (24)$$

with which the source term $f(\boldsymbol{x},t)$ can be derived using Eq. (22).

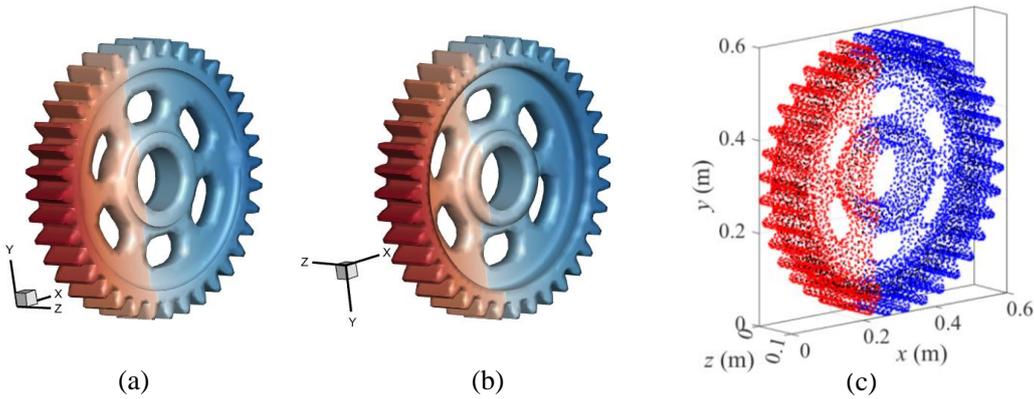

(a)      (b)      (c)

**Fig. 2.** (a) and (b) Geometry of the gear-shaped FGM, and (c) arrangement of collocation nodes (black points) and training nodes (red and blue points) in the SINNs.

In our initial investigation, we examine the viability of the developed SINNs in addressing the dynamic heat conduction problem. For conducting numerical simulations with SINNs, we distribute 2092 collocation nodes (depicted as black points) and 10454 training nodes (indicated by red and blue points) within the gear-shaped domain and on its surface, as illustrated in Fig. 2(c). The architecture of the SINNs comprises two fully-connected hidden layers, each housing 15 neurons. The network undergoes training through 1000 iterations, and the number of Gaussian nodes is taken as $p=5$. Figs. 3(a)-3(c) plot relative errors of the heat fluxes including $u_x(\boldsymbol{x},t)$, $u_y(\boldsymbol{x},t)$ and $u_z(\boldsymbol{x},t)$ on the gear surface at the final time $t=1\text{s}$, attained by the SINNs when Mish function is selected as the activation function. Observations indicate that the heat fluxes predicted by the SINNs exhibit strong agreement with their corresponding exact values, with relative errors concentrated around the orders of -4 and -5. Additionally, Table 1 showcases the performance of the SINNs in solving dynamic problems utilizing different activation functions, such as Sigmoid, Tanh, Swish, Softplus, Arctan and Mish. It is found that all $L_2$ relative errors for the temperatures and heat fluxes are less than $2.49\times10^{-4}$ and $9.04\times10^{-4}$, respectively. These findings demonstrate the feasibility of proposed SINNs for resolving 3D dynamic heat conduction problems in FGM, and also reveal that the performance of the SINNs appears to be relatively insensitive to the choice of activation function.

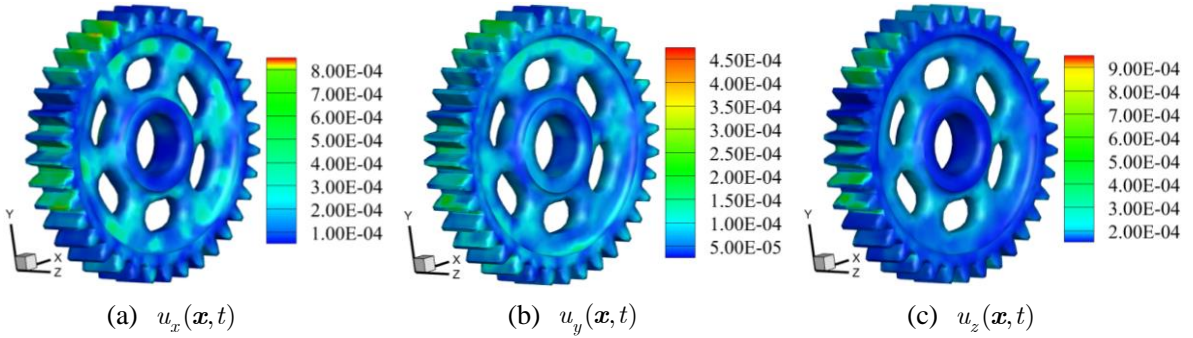

(a) $u_x(\boldsymbol{x},t)$  (b) $u_y(\boldsymbol{x},t)$  (c) $u_z(\boldsymbol{x},t)$

**Fig. 3.** Relative errors of heat fluxes on the gear surface at $t=1\text{s}$ obtained by the SINNs.

**Table 1** Performance of the SINNs using different activation functions.

| Activation function | $u(\boldsymbol{x},t)$ | $u_x(\boldsymbol{x},t)$ | $u_y(\boldsymbol{x},t)$ | $u_z(\boldsymbol{x},t)$ |
|---|---|---|---|---|
| Sigmoid  | $1.82\times10^{-4}$ | $4.72\times10^{-4}$ | $4.96\times10^{-4}$ | $3.39\times10^{-4}$ |
| Tanh     | $2.44\times10^{-4}$ | $6.86\times10^{-4}$ | $6.82\times10^{-4}$ | $9.04\times10^{-4}$ |
| Swish    | $2.82\times10^{-5}$ | $1.50\times10^{-4}$ | $7.55\times10^{-5}$ | $1.10\times10^{-4}$ |
| Softplus | $6.08\times10^{-5}$ | $2.30\times10^{-4}$ | $1.64\times10^{-4}$ | $1.04\times10^{-4}$ |
| Arctan   | $2.49\times10^{-4}$ | $8.24\times10^{-4}$ | $2.73\times10^{-4}$ | $5.27\times10^{-4}$ |
| Mish     | $2.75\times10^{-5}$ | $1.08\times10^{-4}$ | $5.83\times10^{-5}$ | $1.43\times10^{-4}$ |

In the following, we compare the performance of SINNs with that of PINNs in addressing transient heat conduction problems. In both methods, we arrange 1046 collocation nodes and 5227 training nodes within the gear-shaped domain and on its surface. The number of Gaussian nodes is set to $p = 5$ in the SINNs, and similarly, 5 points are evenly distributed over the time interval $(0\text{s}, 1\text{s}]$ in the PINNs. Two fully-connected hidden layers, each consisting of 15 neurons, are employed in both the SINNs and PINNs. The Swish function is used as the activation function in both architectures. Figs. 4(a) and 4(b) show the variation of the loss terms in SINNs and PINNs, respectively, with respect to the number of epoch (The term "epoch" in this paper is different from "iteration", and the maximum value for "epoch" is set to 30000 in all examples). At the initial epoch, the values of the loss terms in both approaches are of the same order of magnitude, while at the final epoch, the values of the loss terms in SINNs and PINNs are approximately $10^{-10}$ and $10^{-8}$, respectively. These results indicate that the SINNs exhibit a faster convergence speed compared to the PINNs. Fig. 5(a) gives the $L_2$ relative errors of temperatures predicted by using the SINNs and the PINNs against the number of iterations, and provides the corresponding CPU-time consumed by both schemes. It is observed that a higher number of iterations leads to improved accuracy in both methods. Additionally, compared to the PINNs, the SINNs achieve higher accuracy with less CPU time, demonstrating significant superior performance. In Fig. 5(b), the relative errors of $u_x(\boldsymbol{x}, t)$ on the gear surface at $t = 1\text{s}$, obtained by utilizing the SINNs and the PINNs through 2500 iterations, reaffirm that the proposed SINNs outperform the PINNs in terms of accuracy.

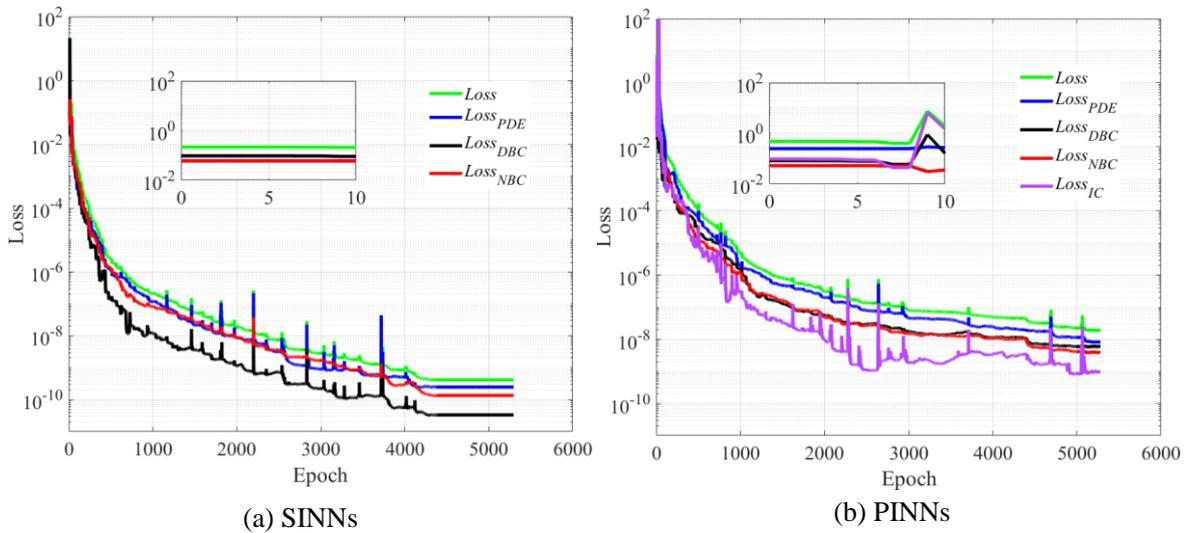

(a) SINNs  (b) PINNs

**Fig. 4.** Loss terms in the PINNs and the SINNs versus number of epoch.

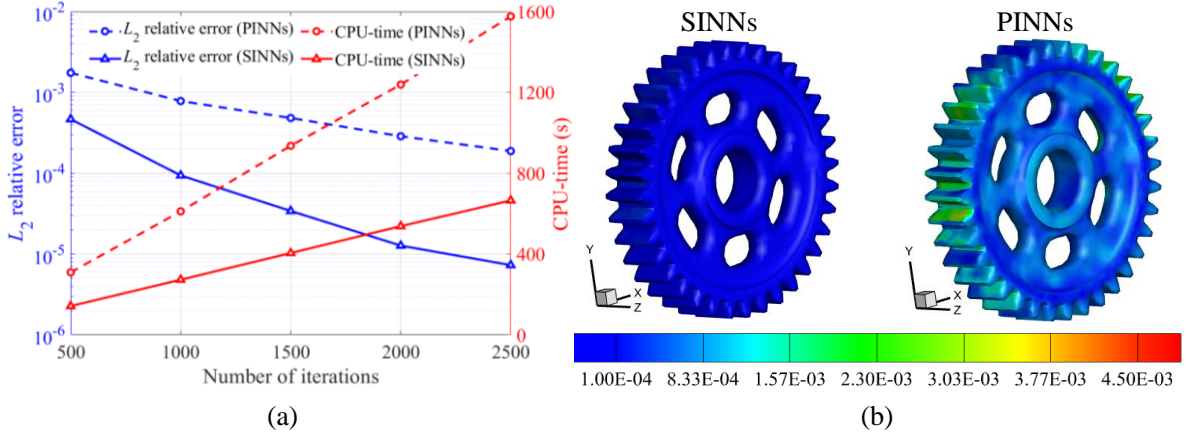

**Fig. 5.** Comparison of the SINNs and the PINNs, (a) $L_2$ relative errors of temperatures and CPU-time with respect to different number of iterations, and (b) relative errors of $u_x(\boldsymbol{x},t)$ on the gear surface at $t=1\text{s}$ using 2500 iterations.

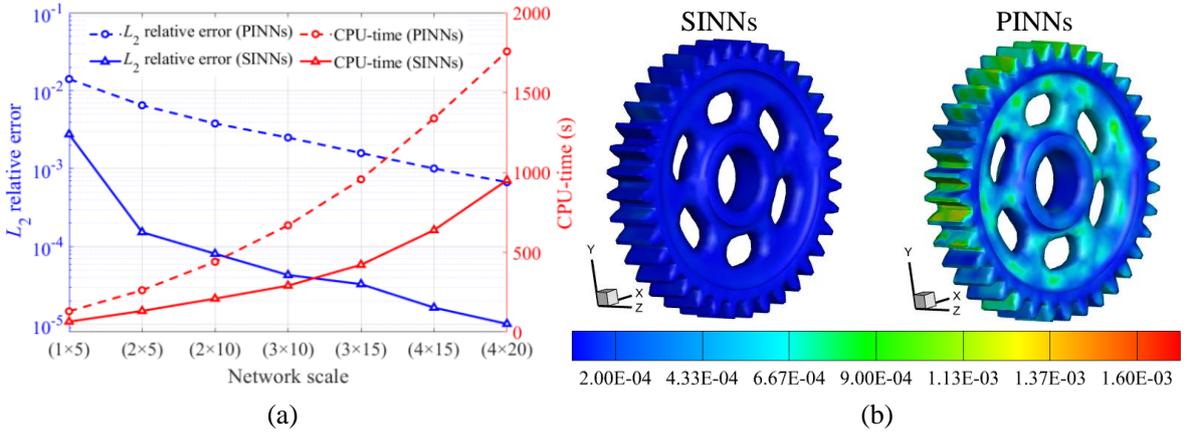

**Fig. 6.** Comparison of the SINNs and the PINNs, (a) $L_2$ relative errors of temperatures and CPU-time in relation to network scales, and (b) relative errors of $u_y(\boldsymbol{x},t)$ on the gear surface at $t=1\text{s}$.

Maintaining the aforementioned node arrangement and parameter settings, we train the SINNs and the PINNs for 1000 iterations, and compare their performance across different network architectures. Fig. 6(a) displays the $L_2$ relative errors of temperatures obtained by employing the SINNs and the PINNs under different network scales, including one fully-connected hidden layer with 5 neurons (1×5), two fully-connected hidden layers with 5 neurons per layer (2×5), two fully-connected hidden layers with 10 neurons per layer (2×10), etc., accompanied by the corresponding CPU-time required for each method. The observation indicates that the SINNs and the PINNs attain better performance as the network scales increase. In addition, for each network architecture, it is noticed that the SINNs consistently achieve higher accuracy compared to the PINNs, and the CPU time required by the SINNs is less than that of the PINNs. Considering the scenario where four fully-connected hidden layers, each consisting of 20 neurons, are employed in both the SINNs and the PINNs, we calculate and present the relative errors of $u_y(\boldsymbol{x},t)$ on the gear surface at $t=1\text{s}$ in Fig. 6(b). The results further affirm the superiority of the developed SINNs over the PINNs in terms of accuracy.

*4.2. Nonlinear heat conduction problem*

In this example, we address a nonlinear heat conduction problem in a mechanical component-shaped domain, as depicted in Figs. 7(a)-7(c). The primary measurements of the mechanical component are $0.50\text{m}$ in length, $0.30\text{m}$ in width, and $0.24\text{m}$ in height. The heat fluxes are specified on the lower surface of the mechanical component ($0\text{m} \leq z \leq 0.12\text{m}$), and the temperatures are stipulated on the remaining surface of the mechanical component. This nonlinear problem is governed by

$$150 u_t(\boldsymbol{x},t) - \nabla \cdot \left[ 0.05 u(\boldsymbol{x},t) + 50 \right] \nabla u(\boldsymbol{x},t) = f(\boldsymbol{x},t), \ \boldsymbol{x} \in \Omega, t \in [0\text{s}, 1\text{s}]. \quad (25)$$

The fabricated solution for this problem is

$$u(\boldsymbol{x},t) = 30\left[\cos(t) + 1.2\right]\left[\sin(x+y) + e^{y+2z}\right], \quad (26)$$

with which the source term $f(\boldsymbol{x},t)$ can be derived using Eq. (25).

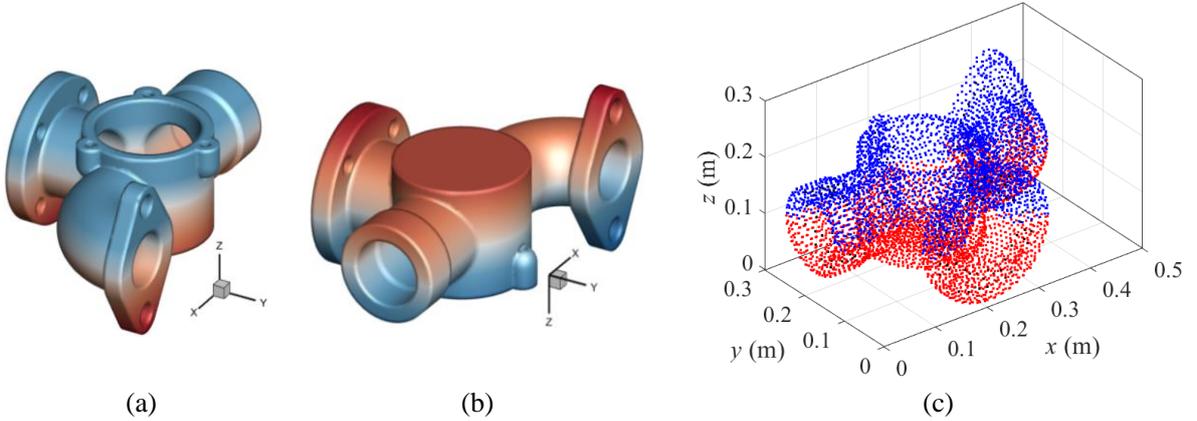

(a) (b) (c)

**Fig. 7.** (a) and (b) Geometry of the mechanical component, and (c) configuration of collocation nodes (black points) and training nodes (red and blue points) in the SINNs.

To address this problem with SINNs, 253 collocation points (depicted as black nodes) and 6322 training points (indicated by red and blue nodes) are strategically arranged inside the mechanical component and along its surface, as illustrated in Fig. 7(c). The architecture of the SINNs consists of two fully-connected hidden layers, each containing 10 neurons. The Swish function is employed as the activation function, and the network undergoes training for 1000 iterations. The number of Gaussian nodes is set to $p = 5$. Fig. 8 presents the comparison results between exact solutions and numerical solutions predicted by utilizing the SINNs. As shown in Fig. 8, the heat fluxes calculated by the SINNs exhibit remarkable consistency with the corresponding exact values, with relative errors primarily concentrated in the orders of -4 and -5. These results demonstrate the feasibility and effectiveness of the proposed SINNs in resolving nonlinear transient heat conduction problems.

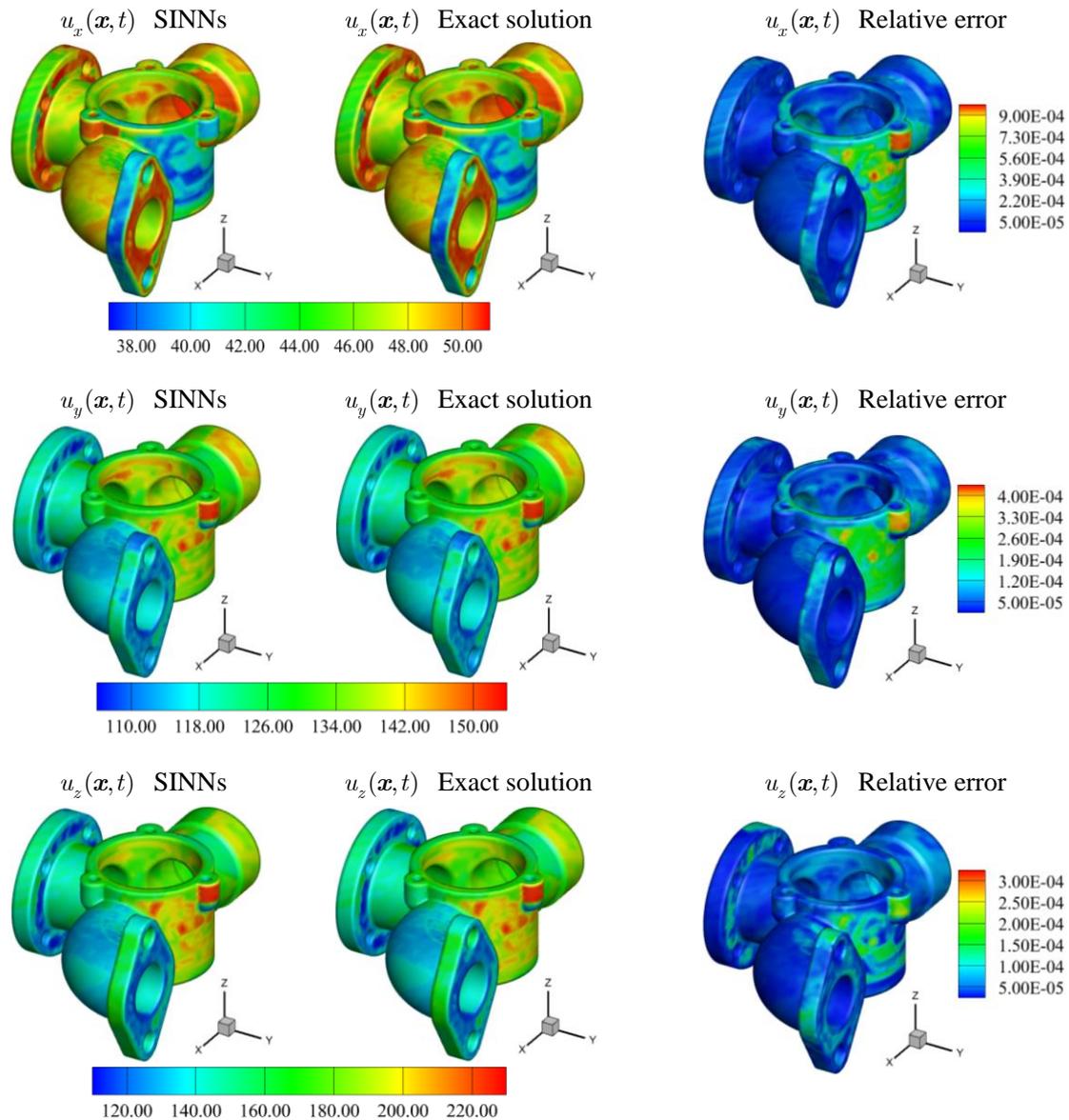

**Fig. 8.** Comparison of exact solutions and numerical solutions obtained by SINNs at $t = 1\text{s}$.

Next, we continue with the following nonlinear heat conduction problem to further examine the performance of the SINNs,

$$150u_t(\boldsymbol{x},t) - \nabla \cdot \left[0.35u(\boldsymbol{x},t) + 20\right]\nabla u(\boldsymbol{x},t) = f(\boldsymbol{x},t), \quad \boldsymbol{x} \in \Omega, t \in [0,T]. \tag{27}$$

The fabricated solution for this problem is

$$u(\boldsymbol{x},t) = 50\left[\cos(t) + 1.2\right]e^{0.3x+0.9y+0.2z}, \tag{28}$$

with which the source term $f(\boldsymbol{x},t)$ can be derived utilizing Eq. (27).

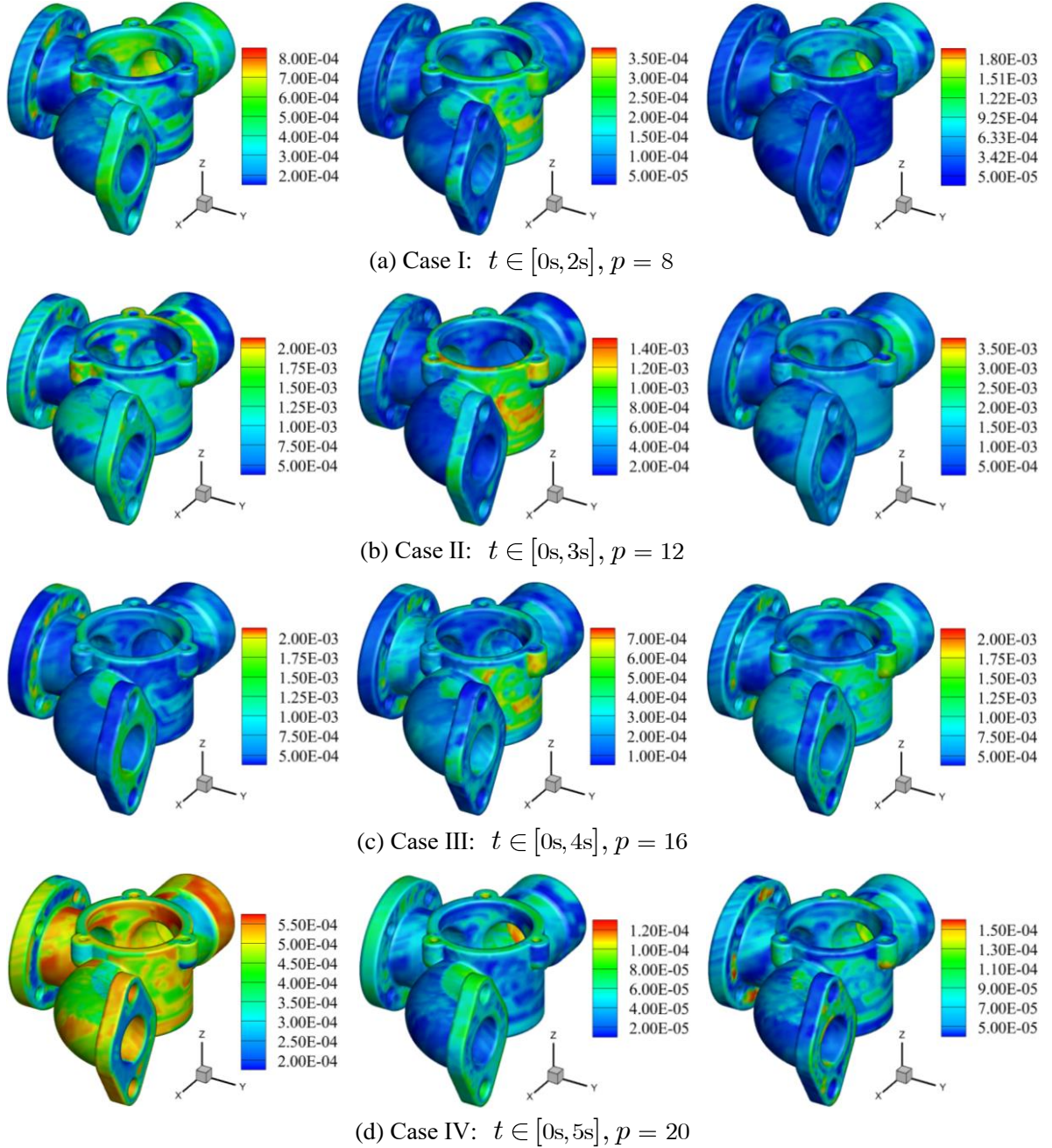

(a) Case I: $t \in [0\text{s}, 2\text{s}]$, $p = 8$

(b) Case II: $t \in [0\text{s}, 3\text{s}]$, $p = 12$

(c) Case III: $t \in [0\text{s}, 4\text{s}]$, $p = 16$

(d) Case IV: $t \in [0\text{s}, 5\text{s}]$, $p = 20$

**Fig. 9.** Relative errors of heat fluxes on the surface of mechanical component for different cases.

In this simulation, we use the same node distribution and parameter settings as employed earlier. Four cases with different time scales are considered, i.e., $[0\text{s}, 2\text{s}]$, $[0\text{s}, 3\text{s}]$, $[0\text{s}, 4\text{s}]$ and $[0\text{s}, 5\text{s}]$. To ensure computational accuracy, the corresponding numbers of Gaussian nodes are assigned as 8, 12, 16, and 20. The relative errors of heat fluxes, including $u_x(\boldsymbol{x},t)$, $u_y(\boldsymbol{x},t)$ and $u_z(\boldsymbol{x},t)$, on the surface of mechanical component for cases I-IV are presented in Figs. 9(a)-9(d), respectively. It is observed that the achieved accuracy is quite satisfactory. The $L_2$ relative errors of temperatures and heat fluxes for different cases are given in Table 2, in which the maximum $L_2$ relative errors for the temperatures and heat fluxes are limited to $2.92 \times 10^{-5}$ and $1.12 \times 10^{-3}$, respectively. The proposed SINNs are proven to be effective for dealing with the nonlinear heat conduction problems, and more attractively, for longer-duration dynamic problems, it can be solved by employing a larger number of Gaussian points.

**Table 2** $L_2$ relative errors of temperatures and heat fluxes for different cases.

| Cases | $u(\boldsymbol{x},t)$ | $u_x(\boldsymbol{x},t)$ | $u_y(\boldsymbol{x},t)$ | $u_z(\boldsymbol{x},t)$ |
|---|---|---|---|---|
| Case I | $1.15 \times 10^{-5}$ | $3.27 \times 10^{-4}$ | $1.33 \times 10^{-4}$ | $5.72 \times 10^{-4}$ |
| Case II | $2.92 \times 10^{-5}$ | $9.35 \times 10^{-4}$ | $5.80 \times 10^{-4}$ | $1.12 \times 10^{-3}$ |
| Case III | $1.56 \times 10^{-5}$ | $6.80 \times 10^{-4}$ | $2.50 \times 10^{-4}$ | $6.76 \times 10^{-4}$ |
| Case IV | $1.54 \times 10^{-5}$ | $4.02 \times 10^{-4}$ | $9.20 \times 10^{-5}$ | $2.74 \times 10^{-4}$ |

*4.3. Wave propagation problem*

In this example, we examine a dynamic wave propagation problem in a cylinder-shaped domain. The diameter of the base of the cylinder is $0.30\text{m}$, and its height is $0.90\text{m}$. Neumann boundary conditions are specified on the top and bottom surfaces of the cylinder, and Dirichlet boundary conditions are stipulated on the remaining surface of the cylinder. The governing equation is described as,

$$u_{tt}(\boldsymbol{x},t) - 250000 \nabla^2 u(\boldsymbol{x},t) = f(\boldsymbol{x},t), \ \boldsymbol{x} \in \Omega, t \in [0\text{s}, 1\text{s}], \tag{29}$$

The fabricated solution for this problem is

$$u(\boldsymbol{x},t) = \sin(2t) e^{x+y+z}, \tag{30}$$

with which the source term $f(\boldsymbol{x},t)$ can be derived utilizing Eq. (29).

We conduct a comparative study on the performance of SINNs and PINNs in handling dynamic wave propagation problems. In both schemes, we arrange 1517 collocation nodes and 2486 training nodes within the cylinder-shaped domain and on its surface. Both architectures utilize three fully-connected hidden layers, each consisting of 10 neurons, with the Mish function employed as the activation function. In the SINNs, the number of Gaussian nodes is taken as $p = 5$, and in the PINNs, 5 points are evenly distributed over the time interval $(0\text{s}, 1\text{s}]$. Fig. 10(a) presents the $L_2$ relative errors in $u(\boldsymbol{x},t)$ over iterations predicted by using SINNs and PINNs, along with the corresponding CPU times for both approaches. It is observed that an increase in the number of iterations enhances the accuracy of both methods. Even more compellingly, the SINNs achieve higher accuracy within a shorter CPU time compared to PINNs, showing a noticeable superiority in performance.

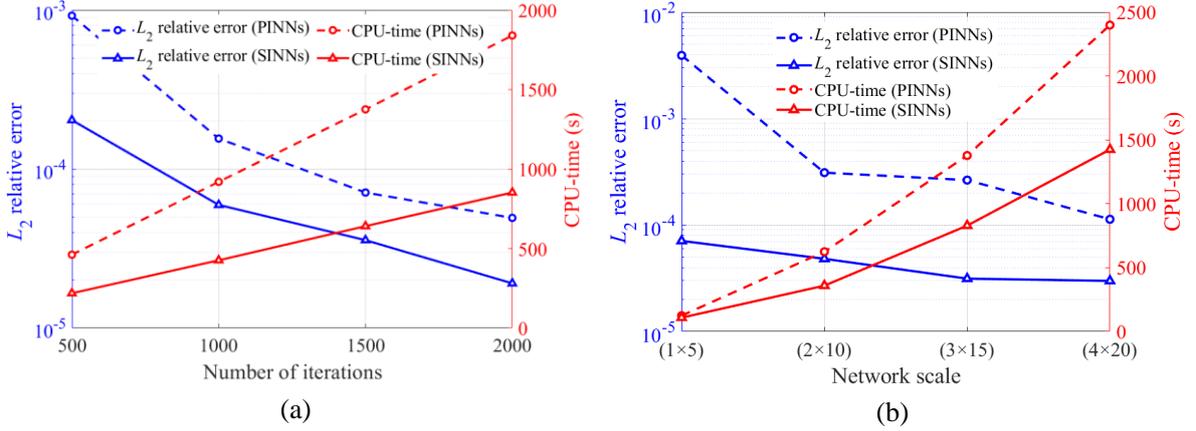

**Fig. 10.** Comparison of the SINNs and the PINNs, $L_2$ relative errors of temperatures and CPU-time with respect to (a) different number of iterations, and (b) different network scales.

Following the previous node distribution and parameter settings, we train the SINNs and PINNs through 1000 iterations and compare their performance in various network architectures. We first consider the scenario where four fully-connected hidden layers, each consisting of 20 neurons, are used in both the SINNs and the PINNs, and respectively depict the variations of the loss terms for SINNs and PINNs with the number of epochs in Figs. 11(a) and 11(b). It is observed that, compared to the PINNs, the SINNs have a faster convergence rate. Next, in Fig. 10(b), we give the $L_2$ relative errors of $u(\boldsymbol{x},t)$ predicted by using SINNs and PINNs under different network scales, including a fully connected hidden layer with 5 neurons (1×5), two fully connected hidden layers with 10 neurons each per layer (2×10), etc. The corresponding CPU time required for each approach is also included. It is found that as the network scale increases, both the SINNs and the PINNs exhibit improved performance. Moreover, in each scenario, the SINNs achieve higher accuracy than the PINNs, with the SINNs also requiring less CPU time. The relative errors of $u_z(\boldsymbol{x},t)$ on the cylinder surface at $t = 1\text{s}$, calculated by the SINNs and the PINNs utilizing different network scales, are illustrated in Fig. 12. The results further validate that the proposed SINNs outperform PINNs in terms of accuracy.

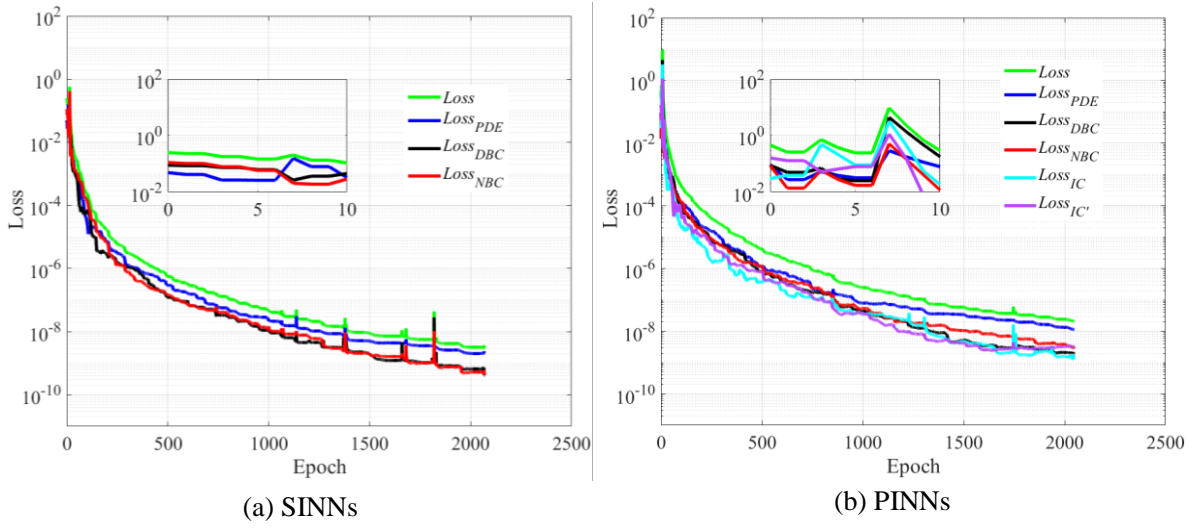

(a) SINNs  (b) PINNs

**Fig. 11.** Loss terms in PINNs and SINNs versus number of epoch.

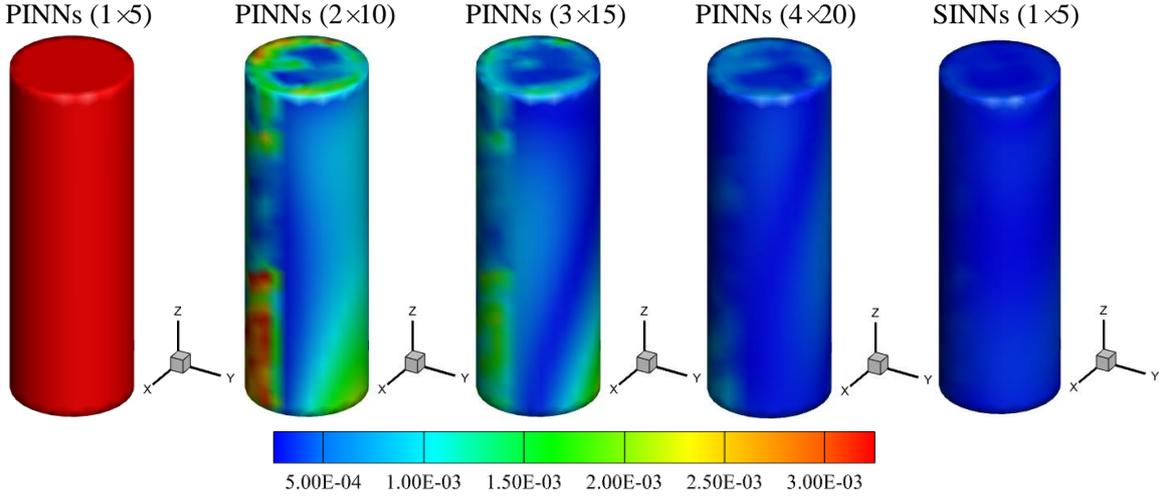

**Fig. 12.** Relative errors of $u_z(\boldsymbol{x},t)$ on the cylinder surface predicted by the SINNs and the PINNs using different network scales.

*4.4. Nonlinear wave propagation problem*

In this example, we address a nonlinear sine-Gordon equation in a dolphin-shaped domain, the primary measurement of the dolphin is $1.84\mathrm{m} \times 0.50\mathrm{m} \times 0.66\mathrm{m}$, as depicted in Figs. 13(a) and 13(b). Dirichlet boundary conditions are specified on the lower surface of the dolphin ($0\mathrm{m} \leq y \leq 0.28\mathrm{m}$), and Neumann boundary conditions are stipulated on the remaining surface of the dolphin. This nonlinear problem is governed by

$$u_{tt}(\boldsymbol{x},t) - \nabla^2 u(\boldsymbol{x},t) + \sin\bigl[u(\boldsymbol{x},t)\bigr] = f(\boldsymbol{x},t),\ \boldsymbol{x} \in \Omega, t \in [0,T], \tag{31}$$

The fabricated solution for this problem is

$$u(\boldsymbol{x},t) = e^{\sin(t)}\left[(2x+y)^2 + e^{y+2z}\right], \tag{32}$$

with which the source term $f(\boldsymbol{x},t)$ can be derived utilizing Eq. (31).

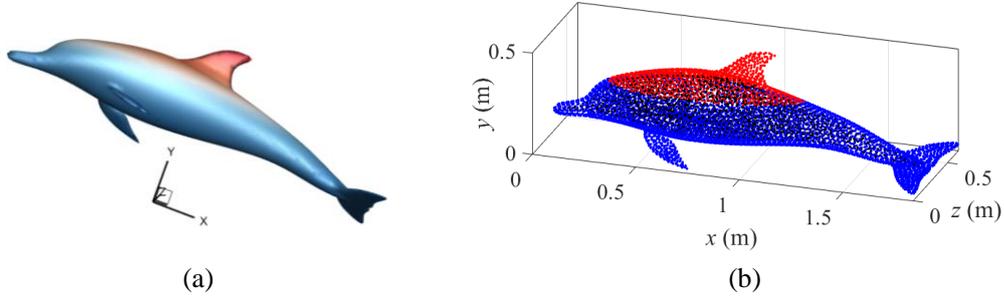

(a)                   (b)

**Fig. 13.** (a) Geometry of the problem, and (b) distribution of collocation nodes (black points) and training nodes (red and blue points) in the SINNs.

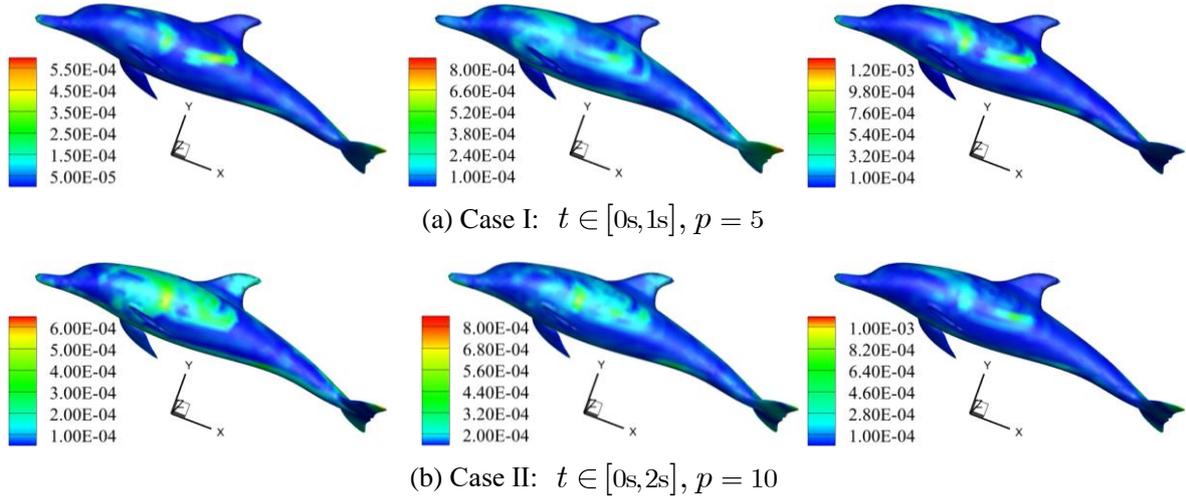

(a) Case I: $t \in [0\text{s}, 1\text{s}]$, $p = 5$

(b) Case II: $t \in [0\text{s}, 2\text{s}]$, $p = 10$

**Fig. 14.** Relative errors of $u_x(\boldsymbol{x},t)$, $u_y(\boldsymbol{x},t)$ and $u_z(\boldsymbol{x},t)$ on the dolphin surface for different cases.

For conducting numerical simulations with SINNs, we distribute 2029 collocation nodes (depicted as black points) and 5278 training nodes (indicated by red and blue points) inside the dolphin-shaped domain and on its surface, as presented in Fig. 13(b). The architecture of the SINNs comprises three fully-connected hidden layers, each housing 15 neurons. The Swish function is employed as the activation function, and the network undergoes training through 1500 iterations. Two cases with different time scales are considered, i.e., $[0\text{s}, 1\text{s}]$ and $[0\text{s}, 2\text{s}]$, the corresponding numbers of Gaussian nodes are set to 5 and 10. Figs. 14(a)-14(b) plot the relative errors of $u_x(\boldsymbol{x},t)$, $u_y(\boldsymbol{x},t)$ and $u_z(\boldsymbol{x},t)$ on the surface of dolphin for cases I and II. The relative errors are found to mainly concentrate in the orders of -4 and -5. The $L_2$ relative errors of $u(\boldsymbol{x},t)$ and its partial derivatives for cases I and II are presented in Table 3, in which all achieved errors are less than $8.52 \times 10^{-5}$. These results demonstrate the effectiveness of the proposed SINNs in addressing nonlinear wave propagation problems.

**Table 3** $L_2$ relative errors of $u(\boldsymbol{x},t)$ and its partial derivatives for different cases.

| Cases | $u(\boldsymbol{x},t)$ | $u_x(\boldsymbol{x},t)$ | $u_y(\boldsymbol{x},t)$ | $u_z(\boldsymbol{x},t)$ |
|---|---|---|---|---|
| Case I | $5.03 \times 10^{-6}$ | $2.16 \times 10^{-5}$ | $6.22 \times 10^{-5}$ | $4.77 \times 10^{-5}$ |
| Case II | $3.22 \times 10^{-5}$ | $8.52 \times 10^{-5}$ | $7.03 \times 10^{-5}$ | $4.65 \times 10^{-5}$ |

*4.5. Inverse problem of heat conduction in FGM*

As the fifth example, we consider an inverse problem of heat conduction in the turbine-shaped FGM, the principal dimension of the turbine is $0.65\text{m} \times 0.31\text{m} \times 0.65\text{m}$, as displayed in Figs. 15(a)-15(c). The governing equation for this inverse problem is described as,

$$36d(\boldsymbol{x})u_t(\boldsymbol{x},t) - \nabla\big[15d(\boldsymbol{x})\nabla u(\boldsymbol{x},t)\big] = f(\boldsymbol{x},t), \ \boldsymbol{x} \in \Omega, t \in [0\text{s},1\text{s}], \tag{33}$$

where the temperature function $u(\boldsymbol{x},t)$ and the function $d(\boldsymbol{x})$ for characterizing material parameters are fabricated as follows,

$$u(\boldsymbol{x},t) = 100\big[\sin(2t) + 1.6\big]e^{0.2x+0.7y+0.1z}, \tag{34}$$

$$d(\boldsymbol{x}) = e^{0.1x} + \sin(y) + z^2. \tag{35}$$

By substituting Eqs. (34) and (35) into Eq. (33), the source term $f(\boldsymbol{x},t)$ can be derived, which will be utilized in the numerical simulation.

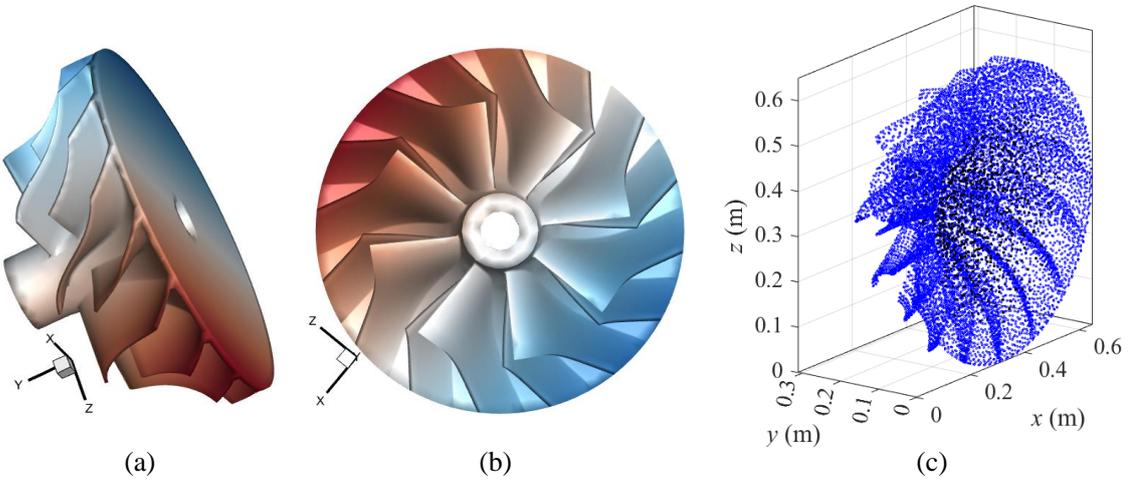

(a) (b) (c)

**Fig. 15.** (a) and (b) Geometry of the turbine-shaped FGM, and (c) arrangement of collocation and training nodes in the SINNs.

To perform numerical simulations using SINNs, we locate 937 collocation nodes (depicted as black points) and 5707 training nodes (indicated by blue points) within the turbine-shaped domain and on its surface, as illustrated in Fig. 15(c). The temperatures are prescribed on the whole surface of the turbine, and some boundary data should be overspecified to solve this inverse problem. In the initial test, the values of heat flux at 20% of the training nodes are known. The architecture of the SINNs consists of three fully-connected hidden layers, each containing 15 neurons. The Swish function is used as the activation function, and the network is trained for 1500 iterations. The number of Gaussian nodes is assigned to $p = 6$. Figs. 16(a) and 16(b) display the correct material parameters, including $\rho(\boldsymbol{x})c(\boldsymbol{x})$ and $\kappa(\boldsymbol{x})$, and the values recovered by the SINNs using different orders of basis functions. Observations show that the predicted values agree well with the corresponding exact ones, and the achieved relative errors are quite satisfactory. Additionally, it can be found that configuring the SINNs with higher-order basis functions results in higher accuracy

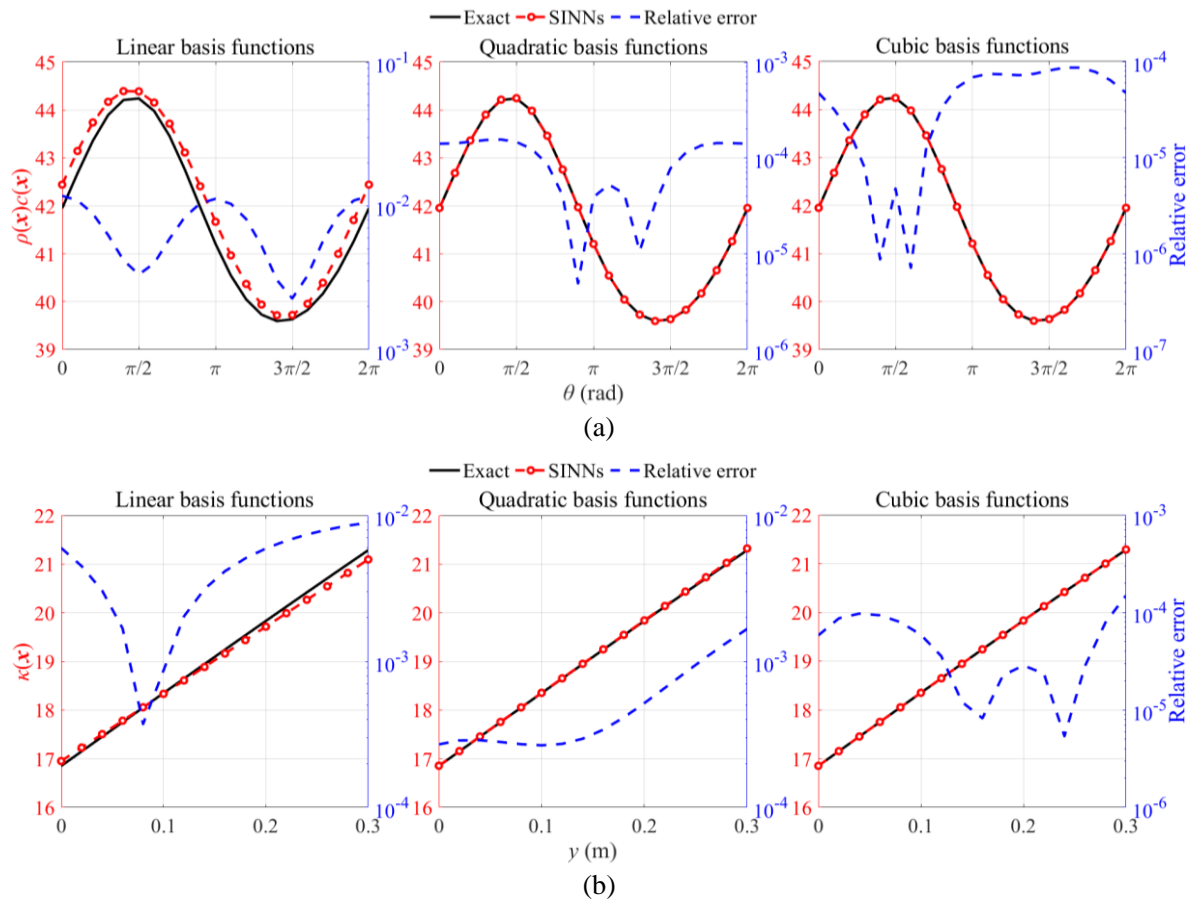

**Fig. 16.** Identified versus correct material parameters (a) $\rho(\boldsymbol{x})c(\boldsymbol{x})$ and (b) $\kappa(\boldsymbol{x})$, and the corresponding relative errors calculated by the SINNs utilizing different orders of basis functions.

Next, we investigate the sensitivity of the developed algorithm with respect to different levels of noise. We use the same node distribution and parameter settings as previously provided, and employ cubic basis functions in this simulation. Fig. 17 depicts the true material parameter $\rho(\boldsymbol{x})c(\boldsymbol{x})$ and the values predicted by the SINNs utilizing boundary data contaminated by noise levels of 1%, 3%, and 5%. The introduction of noise inevitably has an adverse impact on numerical accuracy. However, it is noticeable that the inverted values closely match the true ones. Even with a noise level of 5%, the majority of relative errors still remain below $1.00 \times 10^{-3}$. The relative errors of heat flux $u_x(\boldsymbol{x},t)$ on the turbine surface at $t=1\text{s}$ under different levels of noise are presented in Fig. 18, where the largest error is constrained to $5.00 \times 10^{-3}$. These results confirm that the proposed SINNs can achieve satisfactory accuracy even in the presence of noisy data.

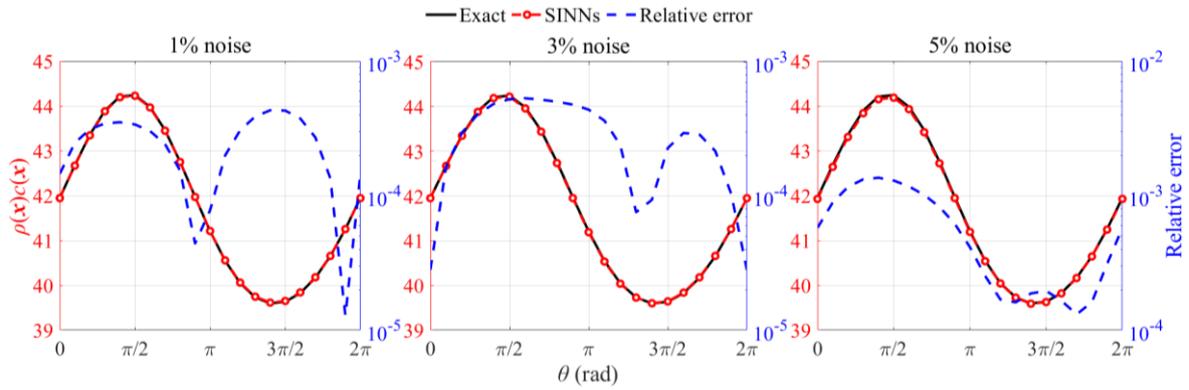

**Fig. 17.** Identified versus correct material parameter $\rho(\boldsymbol{x})c(\boldsymbol{x})$, and the corresponding relative errors calculated by the SINNs under different noise levels.

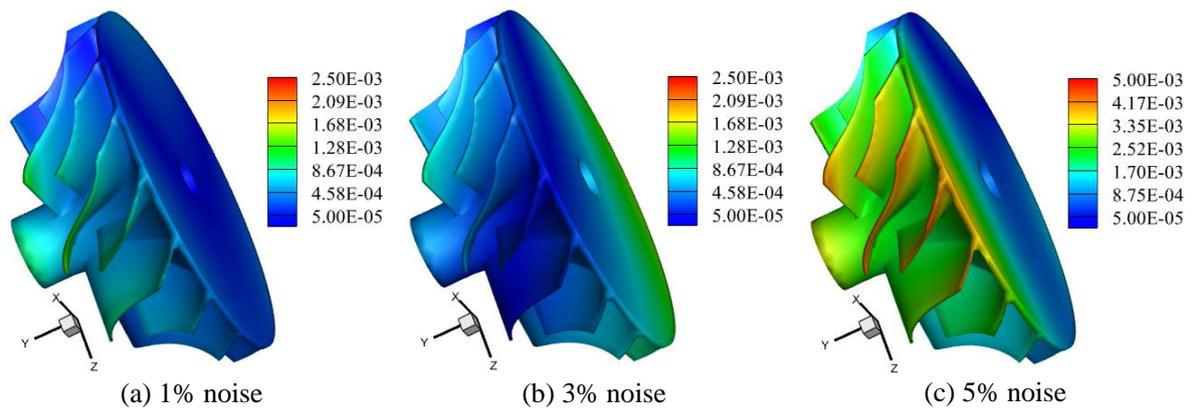

(a) 1% noise     (b) 3% noise     (c) 5% noise

**Fig. 18.** Relative errors of $u_x(\boldsymbol{x},t)$ on the turbine surface at $t=1\text{s}$ under different noise levels.

At the end of this example, we explore the influence of the amounts of overspecified data on the performance of the SINNs. The same node distribution and parameter settings as employed earlier are utilized in this test. In addition, cubic basis functions and boundary data polluted by noise level of 5% are used. Utilizing the SINNs with 15%, 10% and 5% overspecified data, we simultaneously obtain unknown material coefficients and temperature/heat flux distribution. Fig. 19 presents the accurate material parameter $\rho(\boldsymbol{x})c(\boldsymbol{x})$ alongside the values recovered by the SINNs, employing varying amounts of overspecified data. It is found that the predicted values agree well with the true ones in all cases. Fig. 20 displays the relative errors of heat flux $u_x(\boldsymbol{x},t)$ on the turbine surface at $t=1\mathrm{s}$, in which the largest error is limited to $7.00\times10^{-3}$. Notably, the conceived approach can accurately obtain the values of material parameter and heat flux while using a relatively small amount of overspecified data.

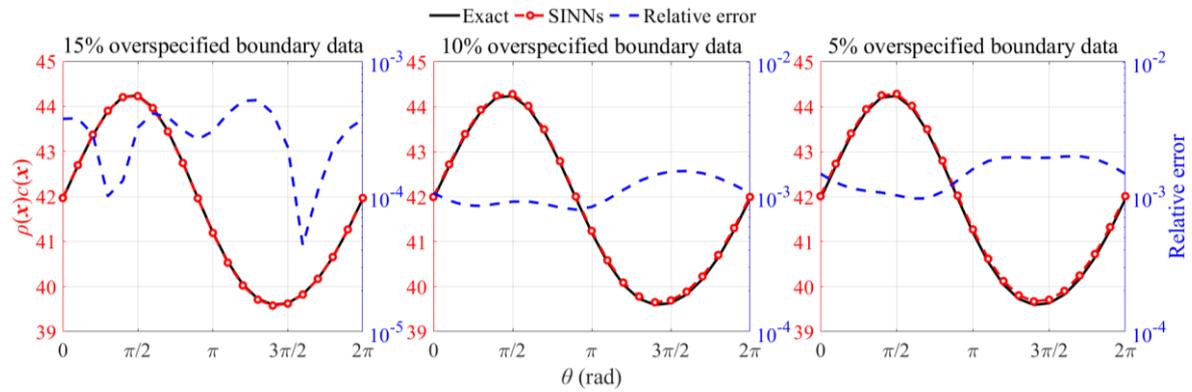

**Fig. 19.** Identified versus correct material parameter $\rho(\boldsymbol{x})c(\boldsymbol{x})$, and the corresponding relative errors obtained by the SINNs using different amounts of overspecified data.

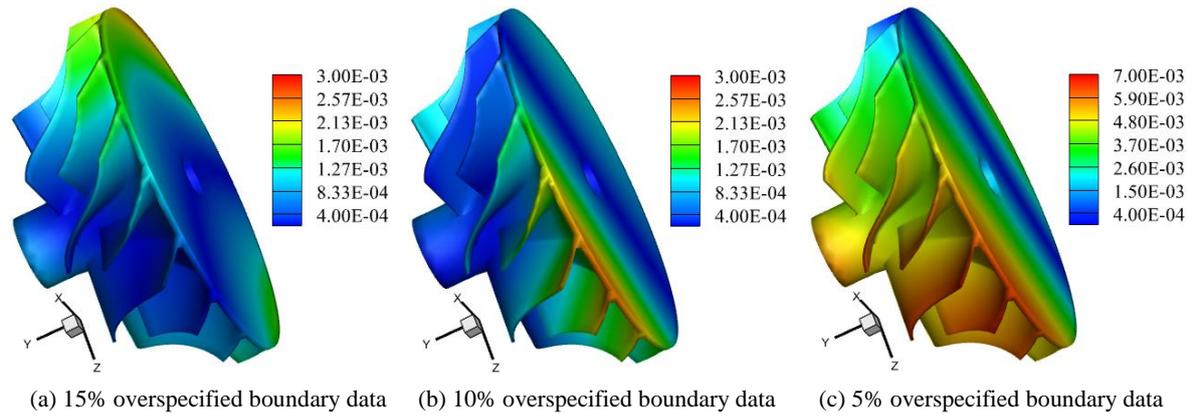

(a) 15% overspecified boundary data    (b) 10% overspecified boundary data    (c) 5% overspecified boundary data

**Fig. 20.** Relative errors of $u_y(\boldsymbol{x},t)$ on the turbine surface at $t=1\mathrm{s}$ utilizing different amounts of overspecified data.

### 4.6. Long-time heat conduction in FGM

Finally, we consider a long-time dynamic problem in the electric motor-shaped FGM, the principal dimension of the electric motor is $0.31\mathrm{m}\times0.63\mathrm{m}\times0.48\mathrm{m}$, as shown in Figs. 21(a)-21(c). The electric motor-shaped FGM is constrained by mixed-type boundary conditions, where the temperatures are enforced on the lower surface of electric motor ( $0\mathrm{m} \leq z \leq 0.31\mathrm{m}$ ), and the heat fluxes are prescribed on the remaining surface of electric motor. The governing equation is the same as Eq. (22) in section 4.1, and the material parameters $\rho(\boldsymbol{x})c(\boldsymbol{x}) = 2.5d(\boldsymbol{x})$, $\kappa(\boldsymbol{x}) = 1.8d(\boldsymbol{x})$, and

$$d(\boldsymbol{x}) = e^{0.6x+0.1y+0.3z}. \tag{36}$$

The fabricated solution for this problem is

$$u(\boldsymbol{x},t) = 25\big[\sin(t) + 1.5\big]e^{0.3x+0.5y+0.2z}, \tag{37}$$

with which the source term $f(\boldsymbol{x},t)$ can be derived utilizing Eq. (22).

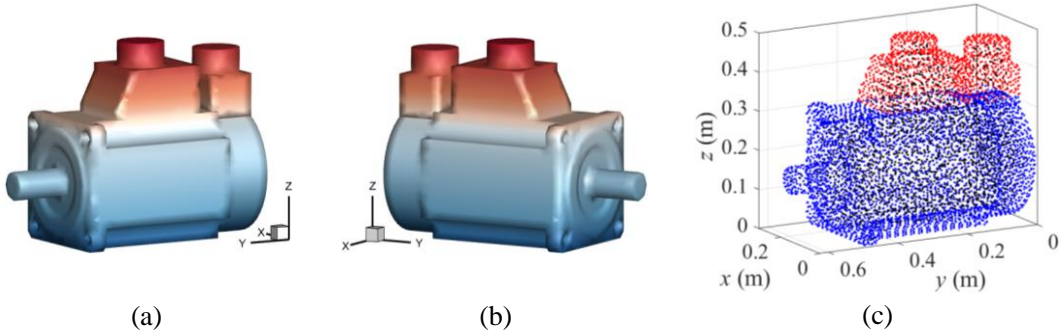

(a)      (b)      (c)

**Fig. 21.** (a) and (b) Geometry of the electric motor-shaped FGM, and (c) arrangement of collocation nodes (black points) and training nodes (red and blue points) in the SINNs.

Using the proposed SINNs, we conduct a long-time simulation from $t = 0\mathrm{s}$ to $t = 100\mathrm{s}$. We first divide this long-time interval into 50 subintervals with a step size of $\Delta T = 2$, and then apply the SINNs to obtain the solution within each subinterval. To solve this heat conduction problem, we distribute 1862 collocation nodes (depicted as black points) and 4850 training nodes (indicated by red and blue points) inside the electric motor-shaped domain and on its surface, as illustrated in Fig. 21(c). The architecture of the SINNs comprises four fully-connected hidden layers, each housing 15 neurons. The Mish function is employed as the activation function, and the network undergoes training for 1500 iterations. The number of Gaussian nodes is taken as $p = 10$ in each subinterval. The $L_2$ relative errors of predicted temperatures and heat fluxes including $u_x(\boldsymbol{x},t)$, $u_y(\boldsymbol{x},t)$ and $u_z(\boldsymbol{x},t)$ at all collocation points are displayed in Fig. 22. It is found that the achieved errors are quite satisfactory, with no remarkable increase observed as the number of time steps increases. Additionally, Figs. 23(a)-23(c) present relative errors of the heat fluxes on the surface of electric motor at the final time $t = 100\mathrm{s}$. Observations indicate that the heat fluxes calculated by the SINNs exhibit strong agreement with their corresponding exact values, with relative errors concentrated around the orders of -4 and -5. These findings suggest that the developed SINNs can effectively and robustly handle long-time dynamic problems.

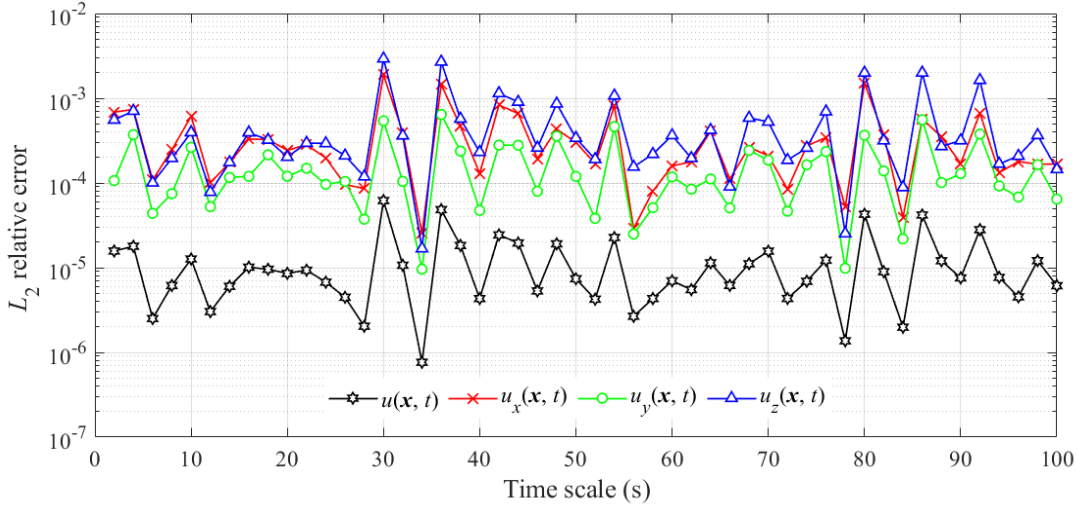

**Fig. 22.** $L_2$ relative error variation of attained temperatures and heat fluxes from $t = 0$s to $t = 100$s.

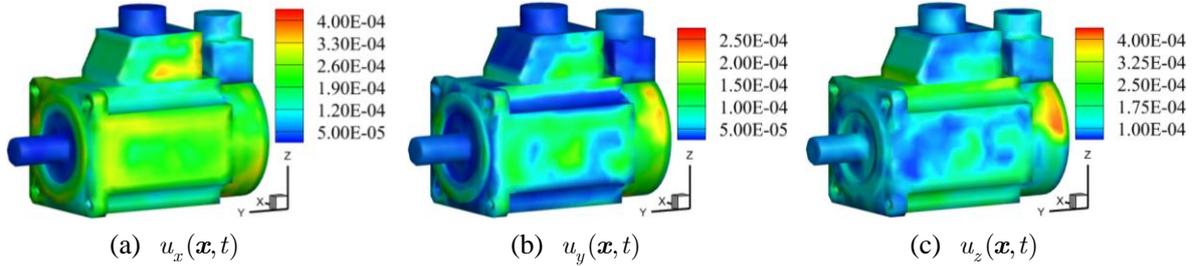

(a) $u_x(\boldsymbol{x},t)$  (b) $u_y(\boldsymbol{x},t)$  (c) $u_z(\boldsymbol{x},t)$

**Fig. 23.** Relative errors of heat fluxes on the surface of electric motor at $t = 100$s calculated by the SINNs.

## 5. Concluding remarks

In this paper, we develop a novel neural network framework, abbreviated as SINNs, to simulate 3D forward and inverse dynamic problems. In the developed methodology, the spectral integration method is adopted for temporal discretization of dynamic problems, and a fully connected neural network with multiple outputs is utilized to approximate solutions in the spatial domain. Using the automatic differentiation technique and spectral integration scheme, we formulate a loss function based on the governing PDEs and boundary conditions. Subsequently, the neural network is trained through the back-propagation of loss function and the gradient descent method. With the view of using the learned network parameters, the SINNs successfully tackle both forward and inverse dynamic problems, encompassing nonlinear PDEs.

Numerical results demonstrate the superior performance of SINNs over the popularly employed PINNs in terms of convergence speed, computational accuracy and efficiency. Additionally, the SINNs exhibit the capability to provide accurate and stable solutions for long-time dynamic problems. It should be noted that this study exclusively examines the performance of the SINNs on some simple and classical dynamic problems. The SINNs also have the potential to tackle more challenging problems, such as those involving thin structures [47], crack propagations [48], and multi-physics coupling [49, 50]. To facilitate the application of the SINNs to real-world issues, some techniques proposed for PINNs, such as domain decomposition [31], self-adaptive method [32], etc., can also be applied to further enhance the accuracy and efficiency of SINNs. Future research endeavors will explore and present results related to these aspects.

**Acknowledgements**

The work described in this paper was supported by the National Natural Science Foundation of China (Grant No. 12302263), the Natural Science Foundation of Qingdao (Grant No. 23-2-1-1-zyyd-jch), the Natural Science Foundation of Shandong Province (Grant Nos. ZR2023QA013, ZR2023YQ005 and ZR2022YQ06), and the Development Plan of Youth Innovation Team in Colleges and Universities of Shandong Province (Grant No. 2022KJ140).

**Declaration of Competing Interest**

The authors declare that they have no known competing financial interests or personal relationships that could have appeared to influence the work reported in this paper.